\documentclass[leqno,12pt]{amsart}
\usepackage{amscd}
\usepackage{latexsym}
\usepackage{epic}

\theoremstyle{plain}
\newtheorem{Theorem}{Theorem}[section]
\newtheorem{Proposition}[Theorem]{Proposition}
\newtheorem{Lemma}[Theorem]{Lemma}
\newtheorem{Corollary}[Theorem]{Corollary}
\newtheorem{Claim}[Theorem]{Claim}

\theoremstyle{definition}

\newtheorem{Remark}[Theorem]{Remark}

\newtheorem{Conjecture}[Theorem]{Conjecture}

\renewcommand{\theTheorem}{\arabic{section}.\arabic{Theorem}}

\newcommand{\ZZ}{{\mathbb{Z}}}
\newcommand{\QQ}{{\mathbb{Q}}}
\newcommand{\RR}{{\mathbb{R}}}

\newcommand{\PP}{{\mathbb{P}}}
\newcommand{\OO}{{\mathcal{O}}}
\newcommand{\rank}{\operatorname{rk}}
\newcommand{\length}{\operatorname{length}}
\newcommand{\ch}{\operatorname{char}}
\newcommand{\dis}{\delta}
\newcommand{\Ker}{\operatorname{Ker}}
\newcommand{\Coker}{\operatorname{Coker}}
\newcommand{\Jac}{\operatorname{Jac}}
\newcommand{\Hilb}{\operatorname{Hilb}}
\newcommand{\Supp}{\operatorname{Supp}}
\newcommand{\Ext}{\operatorname{\mathcal{E}\textsl{xt}}}
\newcommand{\Pic}{\operatorname{Pic}}
\newcommand{\Sing}{\operatorname{Sing}}
\newcommand{\codim}{\operatorname{codim}}
\newcommand{\Bs}{\operatorname{Bs}}
\newcommand{\Res}{\operatorname{Res}}
\newcommand{\Proof}{{\sl Proof.}\quad}
\newcommand{\QED}{{\unskip\nobreak\hfil\penalty50\quad\null\nobreak\hfil
{$\Box$}\parfillskip0pt\finalhyphendemerits0\par\medskip}}
\newcommand{\rest}[2]{\left.{#1}\right\vert_{{#2}}}

\begin{document}

\title[A sharp slope inequality for general stable fibrations of curves]%
{A sharp slope inequality \\
for general stable fibrations of curves}
\author{Atsushi Moriwaki}
\date{1995, December (Version 1.0)}
\address{Department of Mathematics, Faculty of Science,
Kyoto University, Kyoto, 606-01, Japan}
\begin{abstract}
Let $\overline{\mathcal{M}}_g$ be the moduli space of
stable curves of genus $g \geq 2$.
Let $\Delta_i$ be the irreducible component of the boundary of
$\overline{\mathcal{M}}_g$ such that
general points of $\Delta_i$ correspond to stable curves with
one node of type $i$.
Let $\overline{\mathcal{M}}_g^0$ be the set of
stable curves that have
at most one node of type $i>0$.
Let $\delta_i$ be the class of $\Delta_i$ in
$\Pic(\overline{\mathcal{M}}_g) \otimes \QQ$ and $\lambda$ the Hodge class on
$\overline{\mathcal{M}}_g$.
In this paper, we will prove a sharp slope inequality
for general stable fibrations. Namely,
if $C$ is a complete curve on $\overline{\mathcal{M}}_g^0$, then
$\left( (8g+4)\lambda - g \delta_0 -
\sum_{i=1}^{[g/2]} 4i(g-i) \delta_i \ \cdot C \right) \geq 0$.
As an application, we can prove
effective Bogomolov's conjecture for general stable fibrations.
\end{abstract}

\maketitle

\tableofcontents

\section{Introduction}
Let $k$ be an algebraically closed field.
Throughout this paper, we will fix this field $k$.

Let $X$ be a smooth projective surface over $k$,
$Y$ a smooth projective curve over $k$,
and $f : X \to Y$ a generically smooth semistable curve
over $Y$ of genus $g \geq 2$.
Let $P$ be a node of a singular fiber $X_y$ over $y$.
We can assign a number $i$ to the node $P$ in the following way.
Let $\iota : X'_y \to X_y$ be the partial normalization of $X_y$
at $P$. If $X'_y$ is connected, then $i=0$.
Otherwise, $i$ is the minimum of arithmetic genera of
two connected components of $X'_y$. We say the node $P$ of
the singular fiber $X_y$ is of type $i$.
We denote by $\delta_i$ the number of nodes of type $i$ in singular fibers.
If $\ch(k) = 0$,
we know the following inequality due to Cornalba-Harris \cite{CH}
and Xiao \cite{Xi}:
\[
(8g+4)\deg(f_*(\omega_{X/Y})) \geq g \delta,
\]
where $\delta = \delta_0 + \delta_1 + \cdots + \delta_{\left[g/2\right]}$.
(By virtue of \cite{Mo3}, this holds even if $\ch(k) > 0$.)
This inequality is actually sharp because we know
an example which attains equality of the above inequality.
When we consider a fibration with reducible fibers,
we can however observe that
$(8g+4)\deg(f_*(\omega_{X/Y}))$ is rather larger than $g \delta$.
According to the exact formula \cite[Proposition~4.7]{CH} for
hyperelliptic fibrations, we can guess a sharper inequality:
\[
(8g+4)\deg(f_*(\omega_{X/Y})) \geq g \delta_0 +
\sum_{i=1}^{\left[\frac{g}{2}\right]} 4i(g-i) \delta_i.
\]
In this paper, we would like to prove this sharper inequality
for general stable fibrations.

\begin{Theorem}[$\ch(k) = 0$]
\label{thm:sharp:slope:ineq:general:stable:curve}
Let $\bar{f} : \overline{X} \to Y$ be the stable model of
$f : X \to Y$. If every singular fiber of $\bar{f} : \overline{X} \to Y$
has at most one node of type $i > 0$, then
\[
(8g + 4) \deg(f_*(\omega_{X/Y})) \geq
g \delta_0 + \sum_{i=1}^{\left[\frac{g}{2}\right]} 4i(g-i) \delta_i.
\]
In other words, we have the following.
Let $\overline{\mathcal{M}}_g$ be the moduli space of
stable curves of genus $g \geq 2$ over $k$.
Let $\Delta_i$ be the irreducible component of the boundary of
$\overline{\mathcal{M}}_g$ such that
general points of $\Delta_i$ correspond to stable curves with
one node of type $i$.
Let $\overline{\mathcal{M}}_g^0$ be the set of
stable curves that have
at most one node of type $i>0$.
Let $\delta_i$ be the class of $\Delta_i$ in
$\Pic(\overline{\mathcal{M}}_g) \otimes \QQ$
and $\lambda$ the Hodge class on
$\overline{\mathcal{M}}_g$. Then,
for all complete curves $C$ on $\overline{\mathcal{M}}^0_g$,
\[
\left( (8g+4)\lambda - g \delta_0 -
\sum_{i=1}^{\left[\frac{g}{2}\right]} 4i(g-i) \delta_i
\ \cdot \ C \right) \geq 0.
\]
\end{Theorem}

An idea of the proof of
Theorem~\ref{thm:sharp:slope:ineq:general:stable:curve} is as follow.
Let $E$ be the kernel of the natural homomorphism
$f^*f_*(\omega_{X/Y}) \to \omega_{X/Y}$. By \cite{PR},
$E$ is semistable on the geometric generic fiber $X_{\bar{\eta}}$ of $f$.
By the same idea as in \cite{Mo1},
we can apply Bogomolov inequality (Theorem~\ref{thm:B:ineq:fiber:space})
to $E$. This is however
insufficient to get the sharper inequality.
Actually we have only
Cornalba-Harris-Xiao inequality (cf. Remark~\ref{rem:C-H-X:ineq}).
For the sharper inequality, we need to modify $E$ along singular fibers,
namely, we change a compactification of $E_{\bar{\eta}}$ on
$X_{\bar{\eta}}$.
This modification can be done by a special
elementary transformation.

As an application of Theorem~\ref{thm:sharp:slope:ineq:general:stable:curve},
we can show the following answer
for Bogomolov conjecture over function fields
(cf. Theorem~\ref{thm:bogomolov:function:fiels}).
First of all, we fix a notation.
Let $\bar{f} : \overline{X} \to Y$ be the stable model of
$f : X \to Y$. Let $X_y$ (resp. $\overline{X}_y$)
be a singular fiber of $f$ (resp. $\bar{f}$) over $y$, and
$S_y$ the set of nodes $P$ of $\overline{X}_y$ such that
$P$ is not an intersection of two irreducible components of $\overline{X}_y$.
Let $\pi : Z_y \to \overline{X}_y$ be the partial normalization of
$\overline{X}_y$ at each node in $S_y$.
We say $X_y$ is a chain of stable components if
the dual graph of $Z_y$ is homeomorphic to the interval $[0,1]$.

\begin{Theorem}[$\ch(k) = 0$]
We assume that $f$ is not smooth,
every singular fiber of $f$
is a chain of stable components, and one of the following conditions:
\begin{enumerate}
\renewcommand{\labelenumi}{(\alph{enumi})}
\item
the generic fiber of $f$ is hyperelliptic, or

\item
every singular fiber of the stable model $\bar{f} : \overline{X} \to Y$
has at most one node of type $i > 0$.
\end{enumerate}
Then Bogomolov conjecture holds for the generic fiber of $f$,
i.e., we have the following.
Let $K$ be the function field of $Y$, $C$ the generic fiber of $f$,
$\Jac(C)$ the Jacobian of $C$, and let $j : C(\overline{K})
\to \Jac(C)(\overline{K})$ be a morphism given by
$j(P) = \omega_C - P$. Then, $j(C(\overline{K}))$ is
discrete in terms of the semi-norm arising from the Neron-Tate height paring
on $\Jac(C)(\overline{K})$.
More precisely, the maximal radius of ball in which we have only
finitely many points coming from $C(\overline{K})$ via $j$ is
greater than or equal to
\[
\sqrt{\frac{(g-1)^2}{g(2g+1)}\left(
\frac{g-1}{3}\delta_0 + \sum_{i=1}^{\left[\frac{g}{2}\right]}
4i(g-i)\delta_i \right)}.
\]
\end{Theorem}

This is a generalization of the main result of \cite{Mo3} in the
case where $\ch(k) = 0$.
For the proof of the above corollary, we will
calculate an invariant $e(G,D)$ for a metrized graph $G$
with a polarization $D$. Using this,
we have the following inequality, which is rather
weaker than the sharper inequality
(cf. Corollary~\ref{cor:e:for:semistable:chain}).

\begin{Proposition}[$\ch(k) \geq 0$]
If every singular fiber of $f$ is a chain of stable components,
then
\[
\left(\omega_{X/Y} \cdot \omega_{X/Y} \right) \geq
\frac{g-1}{3g} \delta_{0} +
\sum_{i=1}^{\left[\frac{g}{2}\right]}
\left( \frac{4 i (g-i)}{g} - 1 \right) \delta_{i}.
\]
Moreover, if the above inequality is strict,
then Bogomolov conjecture holds for the generic fiber of $f$.
\end{Proposition}

\section{Preliminaries}
\renewcommand{\theTheorem}{\arabic{section}.\arabic{subsection}.\arabic{Theorem}}

In this section, we will recall base loci of canonical linear systems of
reduced curves with only node singularities, and
Bogomolov inequality for semistable vector bundles.

\subsection{Canonical linear systems of
reduced curves with only node singularities}
\setcounter{Theorem}{0}
Let $C$ be a reduced projective curve over $k$ with only node
singularities.
Let $\{ z_1, \ldots, z_s \}$ be the set of singularities of $C$,
and $\pi : \widetilde{C} \to C$ the normalization of $C$.
We set $\pi^{-1}(z_i) = \{ x_i, y_i \}$ for each $i$.
The dualizing sheaf $\omega_C$ of $C$ is defined by
properties: (1) $\omega_C \subset
\pi_*(\Omega^1_{\widetilde{C}}(\sum_i x_i + y_i))$ and
(2) for an open set $U$ of $C$,
$t \in \rest{\omega_C}{U}$ if and only if $\Res_{x_i}(t) + \Res_{y_i}(t) = 0$
for all $i$ with $z_i \in U$.
Note that $\omega_C$ is an invertible sheaf on $C$, i.e.,
a local section $t \in \pi_*(\Omega_{\widetilde{C}}^1(x_i + y_i))_{z_i}$
with $\Res_{x_i}(t) + \Res_{y_i}(t) = 0$ and $\Res_{x_i}(t) \not= 0$
forms a local frame of $(\omega_{C})_{z_i}$.
The arithmetic genus of $C$, denoted by $p_a(C)$, is given by
$\dim_k H^0(C, \omega_C)$ ($=\dim_k H^1(C, \OO_C)$).

\begin{Lemma}[$\ch(k) \geq 0$]
\label{lem:key:for:base:loci}
For smooth points $P_1, \ldots, P_n$ of $C$, let us consider
a homomorphism
\[
  \phi : H^0(C, \omega_C(P_1 + \cdots + P_n)) \to k^n
\]
defined by $\phi(t) = (\Res_{P_1}(t), \ldots, \Res_{P_n}(t))$.
If $C$ is connected, then the kernel of $\phi$ is
$H^0(C, \omega_C)$ and the image of $\phi$ is
the subspace given by
$\{(a_1, \ldots, a_n) \in k^n \mid a_1 + \cdots + a_n = 0 \}$.
\end{Lemma}

\Proof
Let $\pi : \widetilde{C} \to C$, $x_i$, $y_i$ and $z_i$ be
the same as before.
Let $t$ be an element of $H^0(C, \omega_C(P_1 + \cdots + P_n))$
with $\phi(t) = 0$. Then, since $\Res_{P_j}(t) = 0$ for all $j$,
$t$ has no pole at each $P_j$, which implies $t \in H^0(C, \omega_C)$.
Thus, $\Ker(\phi) = H^0(C, \omega_C)$.

Next we would like to see $\sum_{j} \Res_{P_j}(t) = 0$
for all $t \in H^0(C, \omega_C(P_1 + \cdots + P_n))$.
By the definition of $\omega_C$,
$t \in H^0(\widetilde{C},
\Omega^1_{\widetilde{C}}(\sum_i (x_i + y_i) + P_1 + \cdots + P_n))$ and
$\Res_{x_i}(t) + \Res_{y_i}(t) = 0$
for all $i$. On the other hand, by Residue formula,
\[
 \sum_{i=1}^s (\Res_{x_i}(t) + \Res_{y_i}(t)) + \sum_{j=1}^n
 \Res_{P_j}(t) = 0.
\]
Thus, $\sum_{j=1}^n \Res_{P_j}(t) = 0$. Therefore, the image of $\phi$
is contained in $\{(a_1, \ldots, a_n) \in k^n \mid a_1 + \cdots + a_n = 0 \}$.

Since $C$ is connected, by Serre duality,
\[
\dim_k H^1(C, \omega_C(P_1 + \cdots + P_n)) =
\dim_k H^0(C, \OO_C(-P_1 - \cdots -P_n)) = 0.
\]
Thus, using Riemann-Roch theorem, we have
\[
  \dim_k H^0(C, \omega_{C}(P_1 + \cdots + P_n)) =
   \dim_k H^0(C, \omega_{C}) + n - 1.
\]
Since the kernel of $\phi$ is
$H^0(C, \omega_{C})$, the above formula says us that
the dimension of the image of $\phi$ is $n-1$. Therefore,
we get the last assertion.
\QED

\begin{Corollary}[$\ch(k) \geq 0$]
\label{cor:gen:by:glo:irreducible}
Let $C$ be a reduced and irreducible projective curve over $k$ with
only node singularities. If $p_a(C) > 0$, then
$\omega_C$ is generated by global sections.
\end{Corollary}

\Proof
Let $s$ be the number of singularities of $C$ and
$g$ the genus of the normalization of $C$.
If $s=0$ and $g > 0$, or $s=1$ and $g = 0$, then
our assertion is trivial. We will prove it by induction on $s$.
Let $P$ be a node of $C$,
$h : C' \to C$ the partial normalization at $P$, and
$h^{-1}(P) = \{ Q, R \}$. Then, by hypothesis of induction,
$\omega_{C'}$ is generated by global sections and
$H^0(C' \omega_{C'}) \subset H^0(C, \omega_C)$.
Thus, $H^0(C, \omega_{C}) \otimes \OO_C \to \omega_{C}$ is surjective
outside $P$. On the other hand, by Lemma~\ref{lem:key:for:base:loci},
there is a section $t$
of $H^0(C', \omega_{C'}(Q + R))$ such that $\Res_{Q}(t) + \Res_{R}(t) =0$
and $\Res_{Q}(t) \not= 0$. By the definition of $\omega_C$,
$t \in H^0(C, \omega_C)$ and $t$ generates $(\omega_C)_P$.
\QED

\bigskip
Let $C$ be a connected reduced projective curve over $k$
with only node singularities.
Let $P$ be a node of $C$ and $h : C' \to C$ the partial normalization
of $C$ at $P$. We say $P$ is a disconnecting node
if $C'$ is not connected.
Note that, if $C$ is semistable curve, then
a disconnected node is nothing more than a node of type $i > 0$.
An irreducible component $D$ of $C$ is said to be
of socket type if
$D$ is smooth and rational, and all nodes of $C$ on $D$ (i.e.
intersections on $D$ with other components) are disconnecting nodes.
The base locus of $|\omega_C|$, denoted by $\Bs(\omega_C)$,
is defined by the support of
$\Coker(H^0(C, \omega_C) \otimes \OO_C \to \omega_C)$.

\begin{Proposition}[$\ch(k) \geq 0$]
\label{prop:base:locus:stable:curve}
Let $C$ be a connected reduced projective curve over $k$
with only node singularities,
$D_1, \cdots, D_r$ irreducible components of socket type, and
$DN_C$ the set of all disconnecting nodes of $C$. Then we have
the following.
\begin{enumerate}
\renewcommand{\labelenumi}{(\arabic{enumi})}
\item
$\Bs(\omega_C) = D_1 \cup \cdots \cup D_r \cup DN_C$.

\item
Let $P$ be a disconnecting node of $C$,
$h : C' \to C$ the partial normalization of $C$ at $P$,
$C_1$ and $C_2$ two connected components of $C'$.
Then,
\[
H^0(C_1, \omega_{C_1}) \oplus H^0(C_2, \omega_{C_2})
\overset{\sim}{\longrightarrow} H^0(C, \omega_C).
\]

\item
If $P$ is a disconnecting node not lying on
irreducible components of socket type,
then the image of $H^0(C, \omega_C) \otimes \OO_{C,P} \to (\omega_{C})_{P}$
coincides with the image of $(\Omega^1_{C})_{P} \to (\omega_{C})_{P}$.
\end{enumerate}
\end{Proposition}

\Proof
In the following, we denote by $N_C$ the set of all nodes of $C$.
In order to see (1),
it is sufficient to show the following facts:
\begin{enumerate}
\renewcommand{\labelenumi}{(\alph{enumi})}
\item
If $P \in N_C$ and $P$ is a disconnecting node, then
$P \in \Bs(\omega_C)$.

\item
If $P \in N_C$ and $P$ is not a disconnecting node, then
$P \not\in \Bs(\omega_C)$.

\item
If $D$ is an irreducible component of socket type,
then $D \subset \Bs(\omega_C)$.

\item
If $D$ is an irreducible component not of socket type,
then $H^0(C, \omega_C) \otimes \OO_C \to \omega_C$ is surjective
on $D \setminus N_C$.
\end{enumerate}

\medskip
(a) Let $h : C' \to C$ be the partial normalization of $C$ at $P$,
$C_1$ and $C_2$ two connected components of $C'$, and
$h^{-1}(P) = \{ Q, R \}$ with $Q \in C_1$ and $R \in C_2$.
Let $t$ be any section of $H^0(C, \omega_C)$.
Then, $\rest{t}{C_1}$ (resp. $\rest{t}{C_2}$) is a section
of $H^0(C_1, \omega_{C_1}(Q))$ (resp. $H^0(C_2, \omega_{C_2}(R))$).
Then, by Lemma~\ref{lem:key:for:base:loci},
$\Res_{Q}(\rest{t}{C_1}) = \Res_{R}(\rest{t}{C_2}) = 0$.
Thus, $\rest{t}{C_1} \in H^0(C_1, \omega_{C_1})$ and
$\rest{t}{C_2} \in H^0(C_2, \omega_{C_2})$. Hence, $t(P) = 0$.
Therefore, $P \in \Bs(\omega_C)$.

\medskip
(b) Let $h : C' \to C$ be the partial normalization of $C$ at $P$,
and $h^{-1}(P) = \{ Q, R \}$.
By Lemma~\ref{lem:key:for:base:loci}, there is a section $t$
of $H^0(C', \omega_{C'}(Q+R))$ such that $\Res_{Q}(t) + \Res_{R}(t) = 0$
and $\Res_{Q}(t) \not= 0$. Thus, $t$ is a section of
$H^0(C, \omega_C)$ and $t$ generates $(\omega_C)_P$.

\medskip
(c) Let $D \cap N_C = \{ Q_1, \ldots, Q_r \}$.
Let $t$ be any section of $H^0(C, \omega_C)$.
Then, $\rest{t}{D} \in H^0(D, \omega_D(Q_1 + \cdots + Q_r))$.
Since $Q_i$ is a disconnecting node, in the same way as in (a),
we can see that $\Res_{Q_i}(\rest{t}{D}) = 0$ for all $i$.
Thus, $\rest{t}{D}$ has no pole on $D$. Therefore,
$\rest{t}{D} = 0$ because $D \simeq \PP^1$.
Hence, $D \subset \Bs(\omega_C)$.

\medskip
(d) If $p_a(D) > 0$, then our assertion is a consequence of
Corollary~\ref{cor:gen:by:glo:irreducible}.
Thus, we may assume that $D$ is smooth and rational.
Let $E$ be the closure of $C \setminus D$. Since
$D$ is not of socket type, there is a connected component $D'$
of $E$ with $\#(D \cap D') \geq 2$.
We set $D \cap D' = \{ Q_1, \ldots, Q_r \}$ ($r \geq 2$) and
$D'' = E \setminus D'$.
Let $R$ be any point of $D \setminus N_C$. We would like to see
$R \not\in \Bs(\omega_C)$. Since $\omega_D(Q_1 + \cdots + Q_r)$
is generated by global sections, there is a section
$t \in H^0(D, \omega_D(Q_1 + \cdots + Q_r))$ with $t(R) \not= 0$.
Since $\Res_{Q_1}(t) + \cdots + \Res_{Q_r}(t) = 0$,
by Lemma~\ref{lem:key:for:base:loci},
there is a section $t'$ of $H^0(D', \omega_{D'}(Q_1 + \cdots + Q_r))$
such that $\Res_{Q_i}(t') + \Res_{Q_i}(t) = 0$ for all $i$.
Moreover, let $t''$ be the zero form on $D''$.
Then, $t$, $t'$ and $t''$ give a section $s \in H^0(C, \omega_C)$
with $\rest{s}{D} = t$, $\rest{s}{D'} = t'$ and $\rest{s}{D''} = t''$.
Here $s(R) \not= 0$.
Thus, $s$ generates $(\omega_C)_R$.
Therefore, $R \not\in \Bs(\omega_C)$.

\medskip
(2) is obvious by the proof of (a).

\medskip
Finally, let us consider (3).
Let $h : C' \to C$ be the partial normalization of $C$ at $P$,
$C_1$ and $C_2$ two connected components of $C'$, and
$h^{-1}(P) = \{ Q, R \}$ with $Q \in C_1$ and $R \in C_2$.
Let $x$ (resp. $y$) be a local parameter of $C_1$ at $Q$ (resp.
$C_2$ at $R$).
The image of $H^0(C, \omega_C) \otimes \OO_{C,P} \to (\omega_{C})_{P}$
is contained in the image of $(\Omega^1_{C})_{P} \to (\omega_{C})_{P}$
because $H^0(C, \omega_{C}) =
H^0(C_1, \omega_{C_1}) \oplus H^0(C_2,\omega_{C_2})$.
On the other hand,
since $P$ is not lying on irreducible components of socket type,
there are $t_1 \in H^0(C_1, \omega_{C_1})$ and
$t_2 \in H^0(C_2, \omega_{C_2})$ with $t_1(Q) \not= 0$
and $t_2(R) \not= 0$.
Locally, $t_1 = u_1 dx$ and $t_2 = u_2 dy$, where $u_1(Q) \not= 0$ and
$u_2(R) \not= 0$.
Thus, we have (3).
\QED

\subsection{Bogomolov inequality for fiber spaces}
\setcounter{Theorem}{0}
Let $X$ be a smooth projective surface over $k$ and
$E$ a torsion free sheaf of rank $r$ on $X$.
We set
\[
\dis(E) = 2 r c_2(E) - (r -1)c_1(E)^2.
\]
Then we have the following version of Bogomolov inequality
for fiber spaces.

\begin{Theorem}[$\ch(k) = 0$]
\label{thm:B:ineq:fiber:space}
Let $X$ be a smooth projective surface over $k$, $Y$
a smooth projective curve over $k$, and $f : X \to Y$
a surjective morphism with $f_*(\OO_X) = \OO_Y$.
Let $E$ be a torsion free sheaf on $X$.
If $E$ is semistable on the geometric generic fiber $X_{\bar{\eta}}$
of $f$, then $\dis(E) \geq 0$.
\end{Theorem}

\Proof
Let $E^{\vee\vee}$ be the double dual of $E$. Then,
$E^{\vee\vee}$ is locally free, $c_1(E^{\vee\vee}) = c_1(E)$, and
$c_2(E^{\vee\vee}) = c_2(E) - \length(E^{\vee\vee}/E)$.
Thus, we may assume $E$ is locally free.

We assume $\dis(E) < 0$. Then, by Bogomolov instability theorem
(cf. \cite{Bog} or \cite{Mo0}),
there is a non-zero saturated subsheaf $G$ of $E$ such that, if
we set $D = (\rank E)c_1(G) - (\rank G) c_1(E)$, then
$(D^2) > 0$ and $(D \cdot H) > 0$ for some ample divisor $H$.
By Riemann-Roch theorem, we can see that there are
a positive number $n$ and an effective divisor $L$ such that
$nD$ is linearly equivalent to $L$.
Let $F$ be a general fiber of $f$.  Then, $(L \cdot F) \geq 0$, which implies
that $(D \cdot F) \geq 0$. Here we claim that $(D \cdot F) > 0$.
To see this, we assume $(D \cdot F) = 0$. Then, $(L \cdot F) = 0$.
Thus, $L$ is a linear combination of irreducible components of fibers.
By Zariski lemma, this implies $(L^2) \leq 0$, which contradicts to
$(D^2) > 0$. Therefore, we get $(D \cdot F) > 0$.
On the other hand,
\[
(D \cdot F) = (\rank E) \deg(\rest{G}{F}) - (\rank G)\deg(\rest{E}{F}).
\]
Hence, $(D \cdot F) > 0$ means that $\rest{G}{F}$ is a destabilizing subsheaf
of
$\rest{E}{F}$. This contradicts to the assumption that
$E$ is semistable on the geometric generic fiber $X_{\bar{\eta}}$
of $f$.
\QED

\section{Proof of Theorem~\ref{thm:sharp:slope:ineq:general:stable:curve}}
\renewcommand{\theTheorem}{\arabic{section}.\arabic{Theorem}}

To start the proof of Theorem~\ref{thm:sharp:slope:ineq:general:stable:curve},
we need a lot of preparations.

Fix an integer $g \geq 2$ and a polynomial $P_g(n) = (6n-1)(g-1)$.
Let $H_g \subset \Hilb^{P_g}_{\PP^{5g-6}}$ be a subscheme
of all tri-canonically embededded stable curves over $k$,
$Z_g \subset H_g \times \PP^{5g-6}$
the universal tri-canonically embededded stable curves over $k$,
and $\pi : Z_g \to H_g$ the natural projection.
Let $S$ be the set of all points $x \in Z_g$ such that
$\pi$ is not smooth at $x$, and $\Delta = \pi(S)$.
Then, by \cite[Theorem~(1.6) and Corollary~(1.9)]{DM},
$Z_g$, $H_g$, and $S$ are smooth over $k$, $\Delta$ and $\pi^*(\Delta)$
are divisors with only normal crossings, and
$\rest{\pi}{S} : S \to \Delta$ is the normalization of $\Delta$.
Let $\Delta = \Delta_0 \cup \cdots \cup \Delta_{[g/2]}$
be the irreducible decomposition of $\Delta$ such that,
if $x \in \Delta_i \setminus \Sing(\Delta)$, then
$\pi^{-1}(x)$ is a stable curve with one node of type $i$.
Let $\delta_i(x)$ be the number of nodes of type $i$ on a fiber $\pi^{-1}(x)$
over $x \in Z_g$.
Since $\delta_i(x) = \operatorname{mult}_x(\Delta_i)$,
$\delta_i : H_g \to \ZZ$ is upper-semicontinuous.
Therefore, if we set
\[
  H_g^0 = \{ x \in H_g \mid \text{$\pi^{-1}(x)$ has at most one node
of type $i > 0$} \},
\]
then $H_g^0$ is an open set.
In other words, $H_g^0 = H_g \setminus \Sing(\Delta_1 + \cdots
+\Delta_{[g/2]})$.
We set $Z_g^0 = \pi^{-1}(H_g^0)$,
$\Delta^0 = \Delta \cap H_g^0$,
$\Delta_i^0 = \Delta_i \cap H_g^0$,
$S^0 = \left(\rest{\pi}{S}\right)^{-1}(\Delta^0)$,
and $S_i^0 = \left(\rest{\pi}{S}\right)^{-1}(\Delta_i^0)$.
Then, for all $1 \leq i < j \leq [g/2]$,
$\Delta_i^0 \cap \Delta_j^0 = \emptyset$.
For $1 \leq i \leq [g/2]$,
let $\pi^{-1}(\Delta_i^0) = C_i^1 \cup C_i^2$ be
the irreducible decomposition
such that the generic fiber of $\rest{\pi}{C_i^1} : C_i^1 \to \Delta_i^0$
(resp. $\rest{\pi}{C_i^2} : C_i^2 \to \Delta_i^0$) is of genus $i$
(resp. $g-i$).
Then, $C_i^1 \cap C_i^2 = S_i^0$.
Moreover, we denote $S_1^0 \cup \cdots \cup S_{[g/2]}^0$ by
$S^0_{+}$.

\begin{Claim}
\label{claim:image:can:sys}
$\Supp(\Coker(\pi^{*}\pi_{*}(\omega_{Z_g^0/H_g^0}) \to
\omega_{Z_g^0/H_g^0})) = S^0_{+}$ and, for $z \in S^0_{+}$,
the image of $\left(\pi^{*}\pi_{*}(\omega_{Z_g^0/H_g^0})\right)_z \to
\left(\omega_{Z_g^0/H_g^0}\right)_z$ is
$\left(\Omega^1_{Z_g^0/H_g^0}\right)_z$.
\end{Claim}

\Proof
The first assertion is a consequence of (1) of
Proposition~\ref{prop:base:locus:stable:curve}.
Let $F$ be the image of $\pi^{*}\pi_{*}(\omega_{Z_g^0/H_g^0}) \to
\omega_{Z_g^0/H_g^0}$ and $x = \pi(z)$.
As we see in (2) of Proposition~\ref{prop:base:locus:stable:curve},
any section $s$ of $\left(\pi_{*}(\omega_{Z_g^0/H_g^0})\right)_x$
has no pole on irreducible components of fibers around $z$, which means
$F_z \subseteq \left(\Omega^1_{Z_g^0/H_g^0}\right)_z$.
Therefore,
by (3) of Proposition~\ref{prop:base:locus:stable:curve},
we can conclude $F_z = \left(\Omega^1_{Z_g^0/H_g^0}\right)_z$.
\QED

We set
\[
 E = \Ker\left(
\pi^{*}\pi_{*}(\omega_{Z_g^0/H_g^0}) \longrightarrow
\omega_{Z_g^0/H_g^0}\right).
\]

\begin{Claim}
$E$ is locally free.
\end{Claim}

\Proof
Let $I_{S^0}$ (resp. $I_{S^0_{+}}$) be the defining ideal of $S^0$
(resp. $S^0_{+}$). It is well known that
$\Omega^1_{Z_g^0/H_g^0} = I_{S^0} \cdot \omega_{Z_g^0/H_g^0}$
(cf. \cite[Proof of Theorem~5.10]{Mu}).
Thus, by Claim~\ref{claim:image:can:sys},
the image of $\pi^{*}\pi_{*}(\omega_{Z_g^0/H_g^0}) \to
\omega_{Z_g^0/H_g^0}$ is $I_{S^0_{+}} \cdot \omega_{Z_g^0/H_g^0}$.
Since $S^0_{+}$ is smooth and of codimension $2$,
$\Ext^i_{\OO_{Z_g^0}}(I_{S^0_{+}}, \OO_{Z_g^0}) = 0$ for $i \geq 2$.
Therefore, we can see that $\Ext^i_{\OO_{Z_g^0}}(E, \OO_{Z_g^0}) = 0$
for $i \geq 1$. Hence, we have our claim.
\QED

\begin{Claim}
For $1 \leq i \leq [g/2]$,
${\displaystyle
\pi_*(\omega_{C_i^1/\Delta_i^0}) \oplus \pi_*(\omega_{C_i^2/\Delta_i^0})
\overset{\sim}{\longrightarrow}
\pi_*(\omega_{C_i^1 \cup C_i^2/\Delta_i^0}).
}$
\end{Claim}

\Proof
Clearly we have the natural injection
\[
\pi_*(\omega_{C_i^1/\Delta_i^0}) \oplus \pi_*(\omega_{C_i^2/\Delta_i^0})
\longrightarrow \pi_*(\omega_{C_i^1 \cup C_i^2/\Delta_i^0}).
\]
For $x \in \Delta_i^0$, if we set $C_1 = \pi^{-1}(x) \cap C_i^1$ and
$C_2 = \pi^{-1}(x) \cap C_i^2$, then
by (2) of Proposition~\ref{prop:base:locus:stable:curve}
\[
H^0(C_1, \omega_{C_1}) \oplus H^0(C_2, \omega_{C_2}) =
H^0(\pi^{-1}(x), \omega_{\pi^{-1}(x)}).
\]
Thus, the above injection is bijective.
\QED

For $1 \leq i \leq [g/2]$ and $j= 1, 2$, let
\[
  Q_i^j = \Ker\left(
  \left(\rest{\pi}{C_i^j}\right)^*\left(\rest{\pi}{C_i^j}\right)_*
  (\omega_{C_i^j/\Delta_i^0})
  \longrightarrow \omega_{C_i^j/\Delta_i^0}\right).
\]

\begin{Claim}
$Q_i^j$ is a locally free sheaf on $C_i^j$.
\end{Claim}

\Proof
The homomorphism
$\left(\rest{\pi}{C_i^j}\right)^*\left(\rest{\pi}{C_i^j}\right)_*
(\omega_{C_i^j/\Delta_i^0})
\to \omega_{C_i^j/\Delta_i^0}$ is surjective
by (1) of Proposition~\ref{prop:base:locus:stable:curve}.
Therefore, we can see our claim.
\QED

Using projection
$\pi_*(\omega_{C_i^1 \cup C_i^2/\Delta_i^0}) \to
\pi_*(\omega_{C_i^j/\Delta_i^0})$, we have the following commutative diagram:
\[
\begin{CD}
0 @>>> \rest{E}{C_i^j} @>>> \rest{\pi^*\pi_*(\omega_{Z_g^0/H_g^0})}{C_i^j}
  @>>> \rest{\Omega^1_{Z_g^0/H_g^0}}{C_i^j} @>>> 0 \\
@. @VV{\alpha_i^j}V @VVV @VVV @. \\
0 @>>> Q_i^j @>>>
\left(\rest{\pi}{C_i^j}\right)^*\left(\rest{\pi}{C_i^j}\right)_*
(\omega_{C_i^j/\Delta_i^0})
@>>> \omega_{C_i^j/\Delta_i^0} @>>> 0,
\end{CD}
\]
where $\alpha_i^j : \rest{E}{C_i^j} \to Q_i^j$ is the induced
homomorphism.
Thus, we can give a homomorphism
\[
\phi_i : E \to \rest{E}{C_i^1} \oplus \rest{E}{C_i^2}
\overset{\alpha_i^1 \oplus \alpha_i^2}{\longrightarrow} Q_i^1 \oplus Q_i^2.
\]

\begin{Claim}
$\phi_i$ is surjective.
\end{Claim}

\Proof
To show this claim, it is sufficient to see it on each fiber.
Let $C$ be a fiber of $\pi$ over $x \in \Delta_i^0$,
and $C = C^1 \cup C^2$ the irreducible decomposition with
$C^1 = C \cap C_i^1$ and $C^2 = C \cap C_i^2$.
For simplicity, we set $E_C = \rest{E}{C}$,
$Q^1 = \rest{Q_i^1}{C^1}$,
$Q^2 = \rest{Q_i^2}{C^2}$, i.e.,
$E_C$ is the kernel of
$H^0(C, \omega_C) \otimes \OO_C \to \omega_C$ and
$Q^j$ is the kernel of $H^0(C^j, \omega_{C^j}) \otimes \OO_{C^j}
\to \omega_{C^j}$.
Since $H^0(C, \omega_C) = H^0(C^1, \omega_{C^1}) \oplus
H^0(C^2, \omega_{C^2})$,
$Q^1 \oplus Q^2 \subset H^0(C, \omega_{C}) \otimes \OO_C$.
On the other hand, $Q^1, Q^2 \subset E_C$. Thus, we get
$Q^1 \oplus Q^2 \subset E_C$, which shows us that
$E_C \to Q^1 \oplus Q^2$ is surjective.
\QED

Let
\[
F = \Ker\left( \bigoplus_{i=1}^{\left[\frac{g}{2}\right]} \phi_i
\ : \  E \longrightarrow
\bigoplus_{i=1}^{\left[\frac{g}{2}\right]} Q_i^1 \oplus Q_i^2 \right).
\]

\begin{Claim}
$F$ is locally free.
\end{Claim}

\Proof
We set $F_1 = \Ker\left( E \to \bigoplus_{i=1}^{[g/2]} Q_i^1 \right)$.
Then, $F = \Ker\left( F_1 \to \bigoplus_{i=1}^{[g/2]} Q_i^2 \right)$.
$F_1$ is an elementary transformation of $E$.
Thus, by \cite{Ma}, $F_1$ is locally free. Since $\phi_i$ is surjective,
so is $\rest{F_1}{C_i^2} \to Q_i^2$. Therefore,
$F$ is locally free because $F$ is an elementary transformation of $F_1$.
\QED

For a vector bundle $G$ on $Z_g^0$,
we define $\dis(G) \in A_{\dim H_g^0 - 1}(H_g^0)$ by
\[
\dis(G) = \pi_*\left(2 \rank(G) c_2(G) - (\rank(G) - 1)c_1(G)^2\right).
\]

\begin{Claim}
\label{claim:dis:of:E}
${\displaystyle \dis(E) = (8g+4)c_1(\pi_*(\omega_{Z_g^0/H_g^0})) -
g \Delta_0^0 -(3g-2) \sum_{i=1}^{\left[\frac{g}{2}\right]} \Delta_i^0}$.
\end{Claim}

\Proof
First of all, let us recall the Grothendieck-Riemann-Roch theorem.
Let $f : X \to Y$ be a proper morphism of smooth varieties over $k$.
The Grothendieck-Riemann-Roch theorem says that,
for any coherent sheaf $\mathcal{F}$ on $X$,
\[
\operatorname{ch}\left(\sum_{i} (-1)^i R^if_*(\mathcal{F})\right) =
f_*(\operatorname{ch}(\mathcal{F})\operatorname{td}(T_{X/Y})).
\]
We denote $\operatorname{ch}(\sum_{i} R^if_*(\mathcal{F}))$
by $\chi(X/Y, \mathcal{F})$, and $A_i(X)$-component of
$f_*(\operatorname{ch}(\mathcal{F})\operatorname{td}(T_{X/Y}))$
by $f_*(\operatorname{ch}(\mathcal{F})\operatorname{td}(T_{X/Y}))_{(i)}$.
Note that if $\dim(X/Y) = 1$, then
\[
f_*(\operatorname{ch}(\mathcal{F})\operatorname{td}(T_{X/Y}))_{(\dim Y)}
= f_* \left(c_1(\mathcal{F}) -
\rank(\mathcal{F}) \frac{\omega_{X/Y}}{2} \right)
\]
and
\begin{multline*}
f_*(\operatorname{ch}(\mathcal{F})\operatorname{td}(T_{X/Y}))_{(\dim Y - 1)}
= \\
f_* \left(
\frac{c_1(\mathcal{F}) \cdot (c_1(\mathcal{F}) - \omega_{X/Y})}{2}
- c_2(\mathcal{F}) + \rank(\mathcal{F})
\frac{\omega_{X/Y}^2 + c_2(\Omega^1_{X/Y})}{12}
\right).
\end{multline*}

Let us go back to the proof of our claim.
For simplicity, we denote $c_1(\pi_*(\omega_{Z_g^0/H_g^0}))$
and $\omega_{Z_g^0/H_g^0}$
by $\lambda$ and $\omega$ respectively.
Consider an exact sequence
\[
0 \to E \to \pi^*(\pi_*(\omega)) \to
I_{S^0_{+}} \cdot \omega \to 0.
\]
First of all, we have $c_1(E) = \pi^*(\lambda) - \omega$.
Moreover,
\[
\chi(Z_g^0/H_g^0, \pi^*(\pi_*(\omega))) =  \chi(Z_g^0/H_g^0, E) +
\chi(Z_g^0/H_g^0, I_{S^0_{+}} \cdot \omega).
\]
Thus, using the Grothendieck-Riemann-Roch theorem, we can see
\[
\pi_*(c_2(E)) =
\pi_* \left( \omega \cdot \omega - \pi^*(\lambda) \cdot \omega -
S^0_{+} + c_2(\pi^*(\pi_*(\omega))) \right).
\]
Noting that $\pi_*(\pi^*(\lambda) \cdot \omega) = (2g-2)\lambda$,
$\pi_*(c_2(\pi^*(\pi_*(\omega)))) = 0$ and $\pi_*(S_i^0) = \Delta_i^0$,
the above implies
\[
\pi_*(c_2(E)) = \pi_*(\omega \cdot \omega) - (2g-2)\lambda -
\sum_{i=1}^{\left[\frac{g}{2}\right]} \Delta_i^0.
\]
Therefore, by virtue of Noether formula:
$\pi_*(\omega \cdot \omega) = 12 \lambda - \Delta^0$,
we can conclude our claim.
\QED

\begin{Claim}
\label{claim:dis:of:F}
${\displaystyle \dis(F) = (8g+4) c_1(\pi_*(\omega_{Z_g^0/H_g^0})) -
g \Delta_0^0 - \sum_{i=1}^{\left[\frac{g}{2}\right]} 4i(g-i) \Delta_i^0}$.
\end{Claim}

\Proof
Consider an exact sequence:
\[
0 \to F \to E \to \bigoplus_{i=1}^{\left[\frac{g}{2}\right]}
\left( Q_i^1 \oplus Q_i^2 \right)\to 0.
\]
First of all, we have
\begin{align*}
c_1(F) & = c_1(E) -
\sum_{i=1}^{\left[\frac{g}{2}\right]} (\rank(Q_i^1) C_i^1 + \rank(Q_i^2)C_i^2)
\\
& = c_1(E) - \sum_{i=1}^{\left[\frac{g}{2}\right]} ((i-1) C_i^1 + (i'-1)C_i^2),
\end{align*}
where $i' = g-i$.
Moreover, the above exact sequence gives rise to
\[
\chi(Z_g^0/H_g^0, E) = \chi(Z_g^0/H_g^0, F) +
 \sum_{i=1}^{\left[\frac{g}{2}\right]} (\chi(C_i^1/\Delta_i^0, Q_i^1) +
 \chi(C_i^2/\Delta_i^0, Q_i^2)),
\]
which, by Grothendieck-Riemann-Roch theorem, implies
{\allowdisplaybreaks
\begin{multline*}
\pi_*(c_2(F)) = \pi_*\left(c_2(E) +
\frac{c_1(F) \cdot (c_1(F) - \omega)}{2} -
\frac{c_1(E) \cdot (c_1(E) - \omega)}{2}\right) + \\
 \sum_{i=1}^{\left[\frac{g}{2}\right]} \left\{ \left(\rest{\pi}{C_i^1}\right)_*
  \left( c_1(Q_i^1) - \frac{(i-1)\omega_{C_i^1/\Delta_i^0}}{2} \right)
+ \left(\rest{\pi}{C_i^2}\right)_*
  \left( c_1(Q_i^2) - \frac{(i'-1)\omega_{C_i^2/\Delta_i^0}}{2} \right)
\right\}.
\end{multline*}}
Thus, we have
{\allowdisplaybreaks
\begin{multline*}
\dis(F) = \dis(E) + \sum_{i=1}^{\left[\frac{g}{2}\right]}
\pi_*\left(((i-1)C_i^1 + (i'-1)C_i^2)^2\right) \\
+ \sum_{i=1}^{\left[\frac{g}{2}\right]}
\pi_* \left(((g-1) \omega - 2c_1(E))\cdot ((i-1) C_i^1 + (i'-1)C_i^2)\right)
\\
+ 2(g-1) \sum_{i=1}^{\left[\frac{g}{2}\right]} \left(\rest{\pi}{C_i^1}\right)_*
  \left( c_1(Q_i^1) - \frac{(i-1)\omega_{C_i^1/\Delta_i^0}}{2} \right)  \\
+ 2(g-1) \sum_{i=1}^{\left[\frac{g}{2}\right]} \left(\rest{\pi}{C_i^2}\right)_*
  \left( c_1(Q_i^2) - \frac{(i'-1)\omega_{C_i^2/\Delta_i^0}}{2} \right).
\end{multline*}}
Therefore, using formulae:
{\allowdisplaybreaks
\[
\begin{cases}
\left(\rest{\pi}{C_i^1}\right)_*(c_1(Q_i^1)) = -2(i-1)\Delta_i^0, &
\left(\rest{\pi}{C_i^2}\right)_*(c_1(Q_i^2)) = -2(i'-1)\Delta_i^0, \\
& \\
\left(\rest{\pi}{C_i^1}\right)_*(\omega_{C_i^1/\Delta_i^0}) =
2(i-1)\Delta_i^0, &
\left(\rest{\pi}{C_i^2}\right)_*(\omega_{C_i^2/\Delta_i^0}) =
2(i'-1)\Delta_i^0, \\
& \\
\pi_*(c_1(E) \cdot C_i^1) = -(2(i-1) + 1)\Delta_i^0, &
\pi_*(c_1(E) \cdot C_i^2) = -(2(i'-1) + 1)\Delta_i^0, \\
& \\
\pi_*(\omega \cdot C_i^1) = (2(i-1) + 1)\Delta_i^0, &
\pi_*(\omega \cdot C_i^2) = (2(i'-1) + 1)\Delta_i^0, \\
& \\
\pi_*(C_i^1 \cdot C_i^2)  = \Delta_i^0, &
\pi_*(C_i^j \cdot C_i^j)  = -\Delta_i^0 \quad\text{for $j= 1, 2$}, \\
\end{cases}
\]
}
we can see that
\[
\dis(F) = \dis(E) + \sum_{i=1}^{\left[\frac{g}{2}\right]}
(3g-2 - 4ii') \Delta_i^0.
\]
Hence, by Claim~\ref{claim:dis:of:E}, we get
\[
\dis(F) = (8g+4) c_1(\pi_*(\omega_{Z_g^0/H_g^0})) -
g \Delta_0^0 - \sum_{i=1}^{\left[\frac{g}{2}\right]} 4i(g-i) \Delta_i^0.
\]

\textbf{Proof of Theorem~\ref{thm:sharp:slope:ineq:general:stable:curve}:}
\quad
Now we obtain everything to prove
Theorem~\ref{thm:sharp:slope:ineq:general:stable:curve}.
Let $f : X \to Y$ be a semistable curve as in
Theorem~\ref{thm:sharp:slope:ineq:general:stable:curve}.
Then, there is a morphism $h : Y \to H_g^0$ with
$Z_g^0 \times_{H_g^0} Y \simeq \overline{X}$.
Let $h' : X \to Z_g^0$ be the induced morphism.
Let us consider a vector bundle ${h'}^*(F)$ on $Y$.
By \cite{PR}, ${h'}^*(F)$ is semistable on the generic fiber of $f$.
Thus, $\dis({h'}^*(F)) \geq 0$ by Theorem~\ref{thm:B:ineq:fiber:space}.
On the other hand, by Claim~\ref{claim:dis:of:F},
\begin{align*}
\dis({h'}^*(F)) & = \deg(h^*((8g+4) c_1(\pi_*(\omega_{Z_g^0/H_g^0})) -
g \Delta_0^0 - \sum_{i=1}^{\left[\frac{g}{2}\right]} 4i(g-i) \Delta_i^0)) \\
& = (8g+4)\deg(f_*(\omega_{X/Y})) - g \delta_0 -
\sum_{i=1}^{\left[\frac{g}{2}\right]} 4i(g-i) \delta_i.
\end{align*}
Thus, we have
\[
(8g+4)\deg(f_*(\omega_{X/Y})) \geq g \delta_0 +
\sum_{i=1}^{\left[\frac{g}{2}\right]} 4i(g-i) \delta_i.
\]
\QED

\bigskip
Let $\iota : Z_g^0 \to Z_g$ be the inclusion map.
If we set $\overline{F} = \iota_*(F)$, then
$\overline{F}$ is a reflexive coherent sheaf on $Z_g$ because
$\codim(Z_g \setminus Z_g^0) = 2$.
Using $\overline{F}$, we can slightly generalize
Theorem~\ref{thm:sharp:slope:ineq:general:stable:curve}.

\begin{Theorem}[$\ch(k) = 0$]
Let $X$ be a smooth projective surface over $k$,
$Y$ a smooth projective curve over $k$,
and $f : X \to Y$ a semistable curve
of genus $g \geq 2$ over $Y$.
Let $h : Y \to H_g$ and $h' : X \to Z_g$ be the induced morphisms
such that the following diagram is commutative:
\[
\begin{CD}
X @>{h'}>> Z_g \\
@V{f}VV @VV{\pi}V \\
Y @>>{h}> H_g
\end{CD}
\]
If $\overline{F}$ is locally free along $h'(X)$, then
\[
(8g+4)\deg(f_*(\omega_{X/Y})) \geq g \delta_0 +
\sum_{i=1}^{\left[\frac{g}{2}\right]} 4i(g-i) \delta_i.
\]
\end{Theorem}

\Proof
Since $\codim(H_g \setminus H_g^0) = 2$, we have
\[
\dis(\overline{F}) = (8g+4) c_1(\pi_*(\omega_{Z_g/H_g})) -
g \Delta_0 - \sum_{i=1}^{\left[\frac{g}{2}\right]} 4i(g-i) \Delta_i.
\]
Thus, we can conclude our theorem in the same way as before.
\QED

\begin{Remark}
\label{rem:C-H-X:ineq}
Let $X$ be a smooth projective surface over $k$, $Y$
a smooth projective curve over $k$, and $f : X \to Y$
a surjective morphism with $f_*(\OO_X) = \OO_Y$.
Then, in the same idea of Claim~\ref{claim:dis:of:E} or
\cite{Mo1}, we can see that, if $E$ is the kernel of
$f^*f_*(\omega_{X/Y}) \to \omega_{X/Y}$ and
$\omega_{X/Y}$ is $f$-nef,
\[
\delta(E) \leq g(\omega_{X/Y} \cdot \omega_{X/Y})
- 4(g-1)\deg(f_*(\omega_{X/Y})),
\]
where $g \geq 2$ is the genus of the generic fiber.
Thus, by \cite{PR} and Theorem~\ref{thm:B:ineq:fiber:space},
we can recover Cornalba-Harris-Xiao inequality:
\[
g (\omega_{X/Y} \cdot \omega_{X/Y}) \geq
4 (g-1) \deg(f_*(\omega_{X/Y})).
\]
\end{Remark}

\section{Calculation of invariants arising from Green functions}
\label{sec:cal:green:function}
\setlength{\unitlength}{.5in}

In this section, we would like to calculate an invariant
$e(G,D)$ for a metrized graph $G$ with a polarization $D$.
For details of metrized graphs, see \cite{Zh}.

Let $G$ be a connected metrized graph and
$D$ an $\RR$-divisor on $G$.
If $\deg(D) \not= -2$, then
there are a unique measure $\mu_{(G,D)}$ on $G$ and
a unique function $g_{(G,D)}$ on $G \times G$
with the following properties.
\begin{enumerate}
\renewcommand{\labelenumi}{(\alph{enumi})}
\item
${\displaystyle \int_{G} \mu_{(G,D)} = 1}$.

\item
$g_{(G,D)}(x, y)$ is symmetric and continuous on $G \times G$.

\item
For a fixed $x \in G$, $\Delta_y(g_{(G,D)}(x, y)) = \delta_x - \mu_{(G,D)}$.

\item
For a fixed $x \in G$, ${\displaystyle
\int_G g_{(G,D)}(x, y) \mu_{(G,D)}(y) = 0}$.

\item
$g_{(G,D)}(D, y) + g_{(G,D)}(y, y)$ is a constant for all $y \in G$.
\end{enumerate}
The constant $g_{(G,D)}(D, y) + g_{(G,D)}(y, y)$ is denoted by $c(G, D)$.
Further we set
\[
  e(G, D) = 2\deg(D)c(G, D) - g_{(G,D)}(D, D).
\]

First of all, let's consider another expression of $e(G,D)$.

\begin{Lemma}
\label{lem:formula:for:e}
Let $G$ be a connected metrized graph and $D$
an $\RR$-divisor on $G$ with $\deg(D) \not= -2$.
Then, for any point $O$ of $G$,
\[
 e(G,D) = (\deg(D) + 2) g_{(G,D)}(O, D) + r_G(O, D),
\]
where $r_G(P,Q)$ is the resistance between $P$ and $Q$ on $G$.
\end{Lemma}

\Proof
We set $D = \sum_i a_i P_i$. Then,
{\allowdisplaybreaks
\begin{align*}
e(G, D) & = 2 \deg(D) c(G,D) - \sum_i a_i g_{(G,D)}(D, P_i) \\
& = 2\deg(D) c(G,D) - \sum_i a_i \left( c(G,D) - g_{(G,D)}(P_i,P_i) \right) \\
& = \deg(D) c(G,D) + \sum_i a_i g_{(G,D)}(P_i,P_i)
\end{align*}}
Since we know
\[
r_G(P,Q) = g_{(G,D)}(P,P) - 2g_{(G,D)}(P,Q) + g_{(G,D)}(Q,Q)
\]
for all points $P,Q \in G$,
the above implies
{\allowdisplaybreaks
\begin{align*}
e(G,D) & = \deg(D)\left( g_{(G,D)}(O,O) + g_{(G,D)}(O,D) \right) \\
& \qquad + \sum_i a_i \left(
r_G(O,P_i) + 2g_{(G,D)}(O,P_i) - g_{(G,D)}(O,O) \right) \\
& = (\deg(D) + 2)g_{(G,D)}(O,D) + r_G(O,D).
\end{align*}}
\QED

\bigskip
Let $G_1$ and $G_2$ be metrized graphs.
Fix points $x_1 \in G_1$ and $x_2 \in G_2$.
The one point sum $G_1 \vee G_2$ with respect to $x_1$ and $x_2$,
defined by $G_1 \times \{x_2\} \cup \{x_1\} \times G_2$ in $G_1 \times G_2$,
is a metrized graph obtained by joining $x_1 \in G_1$ and $x_2 \in G_2$.
The joining point, which is $\{x_1\}\times\{x_2\}$ in $G_1 \times G_2$,
is denoted by $j(G_1 \vee G_2)$.
Any $\RR$-divisor on $G_i$ ($i=1,2$) can be viewed as an $\RR$-divisor on
$G_1 \vee G_2$.

\begin{Proposition}
\label{prop:e:for:join:graph}
Let $G_1$ and $G_2$ be connected metrized graphs,
and $D_1$ and $D_2$ $\RR$-divisors on $G_1$ and $G_2$ respectively
with $\deg(D_i) \not= -2$ \textup{(}$i=1,2$\textup{)}.
Let $G = G_1 \vee G_2$, $O = j(G_1 \vee G_2)$, and
$D = D_1 + D_2$ on $G_1 \vee G_2$. If $\deg(D_1 + D_2) \not= -2$, then
{\allowdisplaybreaks
\begin{multline*}
e(G, D) = e(G_1, D_1) + e(G_2, D_2) \\
+ \frac{2 \deg(D_2)(\deg(D_1) + 2)g_{(G_1, D_1)}(O,O) +
2 \deg(D_1)(\deg(D_2) + 2)g_{(G_2, D_2)}(O, O)}
{\deg(D_1) + \deg(D_2) + 2}.
\end{multline*}}
Moreover, if $P \in G_2$, then
{\allowdisplaybreaks
\begin{align*}
g_{(G, D)}(P,P) & =
\frac{\deg(D_1)}{\deg(D_1)+\deg(D_2)+2} r_{G_2}(O,P) \\
& \quad +\frac{\deg(D_2) + 2}
{\deg(D_1)+\deg(D_2)+2} g_{(G_2, D_2)}(P,P) \\
& \quad -\frac{\deg(D_1)(\deg(D_2)+2)}
{(\deg(D_1)+\deg(D_2)+2)^2}g_{(G_2,D_2)}(O,O) \\
& \quad +\frac{(\deg(D_1) +2)^2}
{(\deg(D_1) + \deg(D_2) + 2)^2}g_{(G_1, D_1)}(O,O).
\end{align*}}
\end{Proposition}

\Proof
For simplicity, we set $d_i = \deg(D_i)$ and $g_i = g_{(G_i,D_i)}(O,O)$
for $i=1,2$.
By \cite[Lemma~3.7]{Zh}, we have
\[
\mu_{(G,D)} = \frac{d_1+2}{d_1+d_2+2}\mu_{(G_1,D_1)} +
\frac{d_2+2}{d_1+d_2+2}\mu_{(G_2,D_2)} -
\frac{2}{d_1+d_2+2} \delta_O.
\]
Consider the following function on $G$:
\[
g(x) = \begin{cases}
{\displaystyle \frac{d_1+2}{d_1+d_2+2}g_{(G_1,D_1)}(O,x) +
\frac{(d_2+2)^2g_2 - d_2(d_1+2)g_1}{(d_1+d_2+2)^2}}
& \text{if $x \in G_1$}, \\
& \\
{\displaystyle \frac{d_2+2}{d_1+d_2+2}g_{(G_2,D_2)}(O,x) +
\frac{(d_1+2)^2g_1 - d_1(d_2+2)g_2}{(d_1+d_2+2)^2}}
& \text{if $x \in G_2$}. \\
\end{cases}
\]
Then, we can easily check that $g$ is continuous on $G$,
$\Delta(g) = \delta_O - \mu_{(G,D)}$ and $\int_G g \mu_{(G,D)} = 0$.
Thus, $g_{(G,D)}(O,x) = g(x)$. Therefore, by Lemma~\ref{lem:formula:for:e},
we get the first formula. Moreover, using
\[
g_{(G,D)}(P,P) -2g_{(G,D)}(O,P) +
g_{(G,D)}(O,O) = r_G(O,P) = r_{G_2}(O,P)
\]
and
\[
g_{(G_2,D_2)}(P,P) -2g_{(G_2,D_2)}(O,P) +
g_{(G_2,D_2)}(O,O) = r_{G_2}(O,P),
\]
we obtain the second formula.
\QED

\begin{Corollary}
\label{cor:e:for:join:graph:circle}
Let $G$ be a connected metrized graph and $D$ an $\RR$-divisor on $G$ with
$\deg(D) \not= -2$. Let $C$ be a circle of length $l$.
Then,
\[
    e(G \vee C, D) = e(G, D) + \frac{\deg(D)}{3(\deg D + 2)} l.
\]
\end{Corollary}

\Proof
Let $O = j(G \vee C)$ and $t : C \to [0, l)$
a coordinate of $C$ with $t(O) = 0$.
Then, it is easy to see that
\[
  \mu_{(C,0)} = \frac{dt}{l}\quad\text{and}\quad
  g_{(C,0)}(O,x) = \frac{t(x)^2}{2l} - \frac{t(x)}{2} + \frac{l}{12}.
\]
Thus, we have this formula by Proposition~\ref{prop:e:for:join:graph}.
\QED

\bigskip
Next, let's consider $e(G,D)$ for a segment.

\begin{Lemma}
\label{lem:e:for:1:segment}
Let $G$ be a segment of length $l$, and $P$ and $Q$
terminal points of $G$. Let $a$ and $b$ be real numbers
with $a + b \not= 0$, and $D$ an $\RR$-divisor on $G$ given by
$D = (2a-1)P + (2b-1)Q$. Then,
\[
  e(G, D) = \left(\frac{4ab}{a+b} - 1\right)l,\quad
  g_{(G,D)}(P,P) = \frac{b^2}{(a+b)^2}l \quad\text{and}\quad
  g_{(G,D)}(Q,Q) = \frac{a^2}{(a+b)^2}l.
\]
\end{Lemma}

\Proof
First of all, by \cite[Lemma~3.7]{Zh},
\[
 \mu_{(G,D)} = \frac{1}{a+b}(a \delta_P + b \delta_Q).
\]
Let $t : G \to [0, l]$ be a coordinate of $G$.
We set
\[
  f(x) = -\frac{b}{a+b}t(x) + \frac{b^2}{(a+b)^2}l.
\]
Then, $\Delta(f) = \delta_P - \mu$ and $\int_G f\mu = 0$.
Thus, $f(x) = g_{(G,D)}(x, P)$. Therefore,
\[
 g_{(G,D)}(P, P) = \frac{b^2}{(a+b)^2} l
\qquad\text{and}\qquad
 g_{(G,D)}(P, Q) = g_{(G,D)}(Q, P) = -\frac{ab}{(a+b)^2} l.
\]
In the same way, we can see that
\[
  g_{(G,D)}(Q, Q) = \frac{a^2}{(a+b)^2} l.
\]
Therefore, we have
\[
  e(G, D) = \left(\frac{4ab}{a+b} - 1\right)l.
\]
\QED

\bigskip
This lemma can be generalized as follows.

\begin{Proposition}
\label{prop:e:for:n:segment}
Let $G_n$ be a metrized graph given by the following figure.
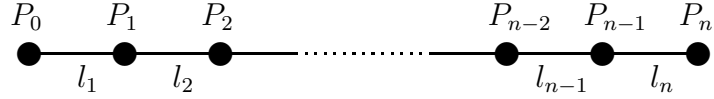
\begin{center}
\begin{picture}(9,2)
\put(1,1){\circle*{.25}}
\put(2,1){\circle*{.25}}
\put(3,1){\circle*{.25}}
\put(6,1){\circle*{.25}}
\put(7,1){\circle*{.25}}
\put(8,1){\circle*{.25}}
\thicklines
\dottedline{0.1}(3.8,1)(5.2,1)
\put(1,1){\line(1,0){1}}
\put(2,1){\line(1,0){1}}
\put(3,1){\line(1,0){0.8}}
\put(6,1){\line(-1,0){0.8}}
\put(6,1){\line(1,0){1}}
\put(7,1){\line(1,0){1}}
\put(0.8,1.3){$P_0$}
\put(1.8,1.3){$P_1$}
\put(2.8,1.3){$P_2$}
\put(5.8,1.3){$P_{n-2}$}
\put(6.8,1.3){$P_{n-1}$}
\put(7.8,1.3){$P_n$}
\put(1.5,0.65){$l_1$}
\put(2.5,0.65){$l_2$}
\put(6.3,0.65){$l_{n-1}$}
\put(7.5,0.65){$l_n$}
\end{picture}
\end{center}
Let $l_i$ be the length between $P_{i-1}$ and $P_i$.
Let
\[
D_n = (2a_0-1)P_0 + (2a_n -1)P_n + \sum_{i=1}^{n-1}2a_iP_i
\]
be an $\RR$-divisor on $G$ with $a_i > 0$ for all $i$. Then, we have
\[
  e(G_n,D_n) = \sum_{i=1}^n \left(
\frac{4(a_0 + \cdots + a_{i-1})(a_{i} + \cdots + a_n)}
{a_0 + \cdots + a_n} - 1 \right)l_i.
\]
\end{Proposition}

\Proof
We set $e_n = e(G_n, D_n)$ and $t_n = g_{(G_n,D_n)}(P_n,P_n)$.
We would like to prove
\[
e_n = \sum_{i=1}^n \left(
\frac{4(a_0 + \cdots + a_{i-1})(a_{i} + \cdots + a_n)}
{a_0 + \cdots + a_n} - 1 \right)l_i
\]
and
\[
t_n = \frac{{\displaystyle
\sum_{i=1}^n \left( a_0 + \cdots + a_{i-1} \right)^2 l_i}}
{(a_0 + \cdots + a_n)^2}.
\]
For this purpose, it is sufficient to show that
\[
t_{n+1} = \frac{(a_0 + \cdots + a_n)^2}{(a_0 + \cdots + a_n + a_{n+1})^2}
(t_n + l_{n+1})
\]
and
\[
e_{n+1} = e_n +
\frac{4a_{n+1}(a_0 + \cdots + a_n)}{a_0 + \cdots + a_n + a_{n+1}}t_n +
\left(
\frac{4 a_{n+1}(a_0 + \cdots + a_n)}{a_0 + \cdots + a_n + a_{n+1}} - 1
\right) l_{n+1}.
\]
Let $L$ be a segment of length $l_{n+1}$, and $Q$ and $P$
terminal points of $L$. Let $E$ be an $\RR$-divisor on $L$ given by
$E = Q + (2a_{n+1}-1)P$.
Let's consider a one point sum $G_n \vee L$ obtained
by joining $P_n$ and $Q$. Then, $G_{n+1} = G_n \vee L$ and
$D_{n+1} = D_n + E$.
Thus, by Proposition~\ref{prop:e:for:join:graph} and
Lemma~\ref{lem:e:for:1:segment},
we have the above recursive equations.
\QED

\begin{Corollary}[$\ch(k) \geq 0$]
\label{cor:e:for:semistable:chain}
Let $X$ be a smooth projective surface over $k$,
$Y$ a smooth projective curve over $k$, and
$f : X \to Y$ a generically smooth semistable curve over $Y$
of genus $g \geq 2$.
Let $X_y$ be a singular fiber of $f$ over $y \in Y$.
Let $X_y = C_1 + \cdots + C_n$ be the irreducible decomposition of $X_y$.
Let $G_y$ be the metrized graph given by the configuration of $X_y$,
$v_i$ the vertex of $G_y$ corresponding to $C_i$, and
$\omega_y$ the divisor on $G_y$ defined by
$\omega_y = \sum_i (\omega_{X/Y} \cdot C_i) v_i$.
If $X_y$ is a chain of stable components, then
\[
e_y = e(G_y, \omega_y) = \frac{g-1}{3g} \delta_{0,y} +
\sum_{i=1}^{\left[\frac{g}{2}\right]}
\left( \frac{4 i (g-i)}{g} - 1 \right) \delta_{i,y},
\]
where $\delta_{i, y}$ is the number of nodes
of type $i$ in $X_y$.
In particular, if every singular fiber of $f$ is
a chain of stable components, then
\[
\left(\omega_{X/Y} \cdot \omega_{X/Y} \right) \geq
\frac{g-1}{3g} \delta_{0} +
\sum_{i=1}^{\left[\frac{g}{2}\right]}
\left( \frac{4 i (g-i)}{g} - 1 \right) \delta_{i}.
\]
Moreover, if the above inequality is strict,
then Bogomolov conjecture holds for the generic fiber of $f$.
\end{Corollary}

\Proof
Under our assumption, $G_y$ can be obtained by
performing one sum of one segment and
many circles. Thus, using Corollary~\ref{cor:e:for:join:graph:circle}
and Proposition~\ref{prop:e:for:n:segment},
we can see the formula for $e_y$.
For the last inequality, note that
$\left(\omega^a_{X/Y} \cdot \omega^a_{X/Y}\right)_a =
\left(\omega_{X/Y} \cdot \omega_{X/Y} \right) - \sum_{y} e_y$
and $\left(\omega^a_{X/Y} \cdot \omega^a_{X/Y}\right)_a \geq 0$
(cf. \cite{Zh} or \cite{Mo3}).
Furthermore,
if $\left(\omega^a_{X/Y} \cdot \omega^a_{X/Y}\right)_a > 0$,
then Bogomolov conjecture holds for the generic fiber of $f$. (cf. \cite{Zh})
\QED

\section{Bogomolov conjecture over function fields}

Let $X$ be a smooth projective surface over $k$,
$Y$ a smooth projective curve over $k$,
and $f : X \to Y$ a generically smooth semistable curve
of genus $g \geq 2$ over $Y$.
Let $K$ be the function field of $Y$, $\overline{K}$
the algebraic closure of $K$, and $C$ the generic fiber of $f$.
Let $j : C(\overline{K}) \to \Jac(C)(\overline{K})$ be a morphism
given by $j(x) = (2g-2)x - \omega_C$ and $\Vert\ \Vert_{NT}$
the semi-norm arising from the Neron-Tate height pairing on
$\Jac(C)(\overline{K})$.
We set
\[
B_C(P;r) = \left\{ x \in C(\overline{K}) \mid
\Vert j(x) - P \Vert_{NT} \leq r \right\}
\]
for $P \in \Jac(C)(\overline{K})$ and $r \geq 0$, and
\[
r_C(P) =
\begin{cases}
-\infty & \mbox{if $\#\left(B_C(P;0)\right) = \infty$}, \\
& \\
\sup \left\{ r \geq 0 \mid \#\left(B_C(P;r)\right) < \infty \right\} &
\mbox{otherwise}.
\end{cases}
\]
An effective version of Bogomolov conjecture claims the following.

\begin{Conjecture}[Effective Bogomolov conjecture]
\label{conj:effective:bogomolov}
If $f$ is non-isotrivial, then
there is an effectively calculated positive number
$r_0$ with
\[
\inf_{P \in \Jac(C)(\overline{K})} r_C(P) \geq r_0.
\]
\end{Conjecture}

In \cite{Mo3} and \cite{Mo4}, we proved the following results.

\begin{enumerate}
\renewcommand{\labelenumi}{(\Alph{enumi})}
\item ($\ch(k) \geq 0$)
If $f$ is non-isotrivial and the stable model of $f : X \to Y$
has only irreducible fibers, then
\[
\inf_{P \in \Jac(C)(\overline{K})} r_C(P) \geq
\begin{cases}
\sqrt{12(g-1)} & \mbox{if $f$ is smooth}, \\
& \\
{\displaystyle \sqrt{\frac{(g-1)^3}{3g(2g+1)}\delta_0}} & \mbox{otherwise}.
\end{cases}
\]

\item ($\ch(k) \geq 0$)
If $f$ is non-isotrivial and $g = 2$, then
$f$ is not smooth and
\[
\inf_{P \in \Jac(C)(\overline{K})} r_C(P) \geq
\sqrt{\frac{2}{135} \delta_0 + \frac{2}{5} \delta_1}.
\]
(According to the exact calculations in \cite{Mo4},
$(\omega_{X/Y} \cdot \omega_{X/Y}) =
\frac{1}{5}\delta_0 + \frac{7}{5} \delta_1$ and
$\sum_y e_y \leq \frac{5}{27}\delta_0 + \delta_1$.)
\end{enumerate}

In this section, we would like to prove the following answer
as an application of our slope inequality.

\begin{Theorem}[$\ch(k) = 0$]
\label{thm:bogomolov:function:fiels}
We assume that $f$ is not smooth,
every singular fiber of $f$
is a chain of stable components, and one of the following conditions:
\begin{enumerate}
\renewcommand{\labelenumi}{(\alph{enumi})}
\item
the generic fiber of $f$ is hyperelliptic, or

\item
every singular fiber of the stable model of $f : X \to Y$
has at most one node of type $i>0$.
\end{enumerate}
Then we have
\[
\inf_{P \in \Jac(C)(\overline{K})} r_C(P) \geq
\sqrt{\frac{(g-1)^2}{g(2g+1)}\left(
\frac{g-1}{3}\delta_0 + \sum_{i=1}^{\left[\frac{g}{2}\right]}
4i(g-i)\delta_i \right)}.
\]
\end{Theorem}

\Proof
First of all, note the following fact
(cf. \cite[Theorem 5.6]{Zh},
\cite[orollary 2.3]{Mo3} or \cite[Theorem 2.1]{Mo4}).
If $(\omega_{X/Y}^a \cdot \omega_{X/Y}^a)_a > 0$, then
\[
\inf_{P \in \Jac(C)(\overline{K})} r_C(P) \geq
\sqrt{(g-1)(\omega_{X/Y}^a \cdot \omega_{X/Y}^a)_a},
\]
where $(\ \cdot \ )_a$ is the admissible pairing.

By the definition of admissible pairing,
we can set
\[
\left(\omega_{X/Y}^a \cdot \omega_{X/Y}^a \right)_a =
\left(\omega_{X/Y} \cdot \omega_{X/Y} \right) - \sum_{y \in Y} e_y,
\]
where $e_y$ is $e(G_y, \omega_y)$
treated in \S\ref{sec:cal:green:function}.
Under our assumption, by \cite[Proposition~4.7]{CH} and
Theorem~\ref{thm:sharp:slope:ineq:general:stable:curve},
we have
\[
(8g + 4) \deg(f_*(\omega_{X/Y})) \geq
g \delta_0 + \sum_{i=1}^{\left[\frac{g}{2}\right]} 4i(g-i) \delta_i.
\]
Thus, using Noether formula, the above inequality implies
\[
(\omega_{X/Y} \cdot \omega_{X/Y}) \geq
\frac{g-1}{2g+1} \delta_0 + \sum_{i=1}^{\left[\frac{g}{2}\right]}
\left(\frac{12i(g-i)}{2g+1} - 1 \right) \delta_i.
\]
Moreover,
by Corollary~\ref{cor:e:for:semistable:chain}, we get
\[
\sum_{y} e_y = \frac{g-1}{3g} \delta_{0} +
\sum_{i=1}^{\left[\frac{g}{2}\right]}
\left( \frac{4i(g-i)}{g} - 1 \right) \delta_{i}.
\]
Thus, we have our theorem.
\QED


\end{document}